\pdfoutput=1
\documentclass{ifacconf}

\usepackage{graphicx}      
\usepackage{natbib}        

\usepackage{xcolor}
\usepackage{pgfplots}
\pgfplotsset{compat=1.18}

\definecolor{viridisYellow}{RGB}{253,231,37}
\definecolor{viridisGreen}{RGB}{94,201,98}
\definecolor{viridisTeal}{RGB}{33,145,140}
\definecolor{viridisBlue}{RGB}{59,82,139}
\definecolor{viridisViolet}{RGB}{68,1,84}
\definecolor{matchingRed}{HTML}{D01F3C}
\definecolor{matchingOrange}{HTML}{FFC077}

\usetikzlibrary{shapes.geometric, arrows, positioning, calc, decorations.pathreplacing, tikzmark, patterns, arrows.meta, shapes.symbols}
\pgfdeclarelayer{background}
\pgfdeclarelayer{foreground}
\pgfsetlayers{background,main,foreground}

\tikzstyle{model} = [rounded corners=0.25cm, minimum height=.5cm, text centered, draw=matchingOrange, fill=matchingOrange!50, font=\small, minimum width=2cm]
\tikzstyle{error-estimator} = [diamond, text centered, draw=viridisGreen, fill=viridisGreen!50, aspect=1.5, inner sep=1pt, font=\small]
\tikzstyle{output} = [rectangle, minimum width=2.2cm, minimum height=.5cm, text centered, draw=viridisYellow, fill=viridisYellow!50, font=\small]
\tikzstyle{arrow} = [thick, ->, >=stealth]

\newcommand{\distbg}{0.15}
\colorlet{colorbg}{viridisBlue}

\newcommand{\inputSpace}{\mathcal{P}}
\newcommand{\complexity}[1]{\mathcal{C}(#1)}
\newcommand{\errorMeasure}[1]{\mathcal{E}(#1)}
\newcommand{\solutionMap}{S}
\newcommand{\solutionSpace}{\mathcal{V}}

\begin{document}
\begin{frontmatter}

\title{Adaptive Model Hierarchies\\for Multi-Query Scenarios\thanksref{footnoteinfo}} 

\thanks[footnoteinfo]{The authors acknowledge funding by the Deutsche Forschungsgemeinschaft (DFG, German Research Foundation) under Germany's Excellence Strategy EXC 2044 –390685587, Mathematics Münster: Dynamics–Geometry–Structure.}

\author[First]{Hendrik Kleikamp}
\author[First]{Mario Ohlberger}

\address[First]{Institute for Analysis and Numerics, Mathematics Münster, University of Münster, Einsteinstrasse 62, 48149 Münster, Germany (e-mail: hendrik.kleikamp@uni-muenster.de).}

\end{frontmatter}

\section{Introduction}
In this contribution we present an abstract framework for adaptive model hierarchies together with several instances of hierarchies for specific applications.
The hierarchy is particularly useful when integrated within an outer loop, for instance an optimization iteration or a Monte Carlo estimation where for a large set of requests answers fulfilling certain criteria are required.
Within the hierarchy, multiple models are combined and interact with each other pursuing the overall goal to reduce the run time in a multi-query context.
To this end, models with different accuracies and effort for evaluation are used in such a way that the cheapest (and typically least accurate) models are evaluated first when a request comes in.
If the result fulfills a prescribed criterion, it can be returned to the outer loop.
Otherwise, the model hierarchy falls back to more costly, but at the same time more accurate, models.
The cheaper models are improved by means of training data gather whenever the more accurate models are evaluated.
\par
In the next section we provide an abstract and detailed description of the components of the hierarchy and their interaction.
Subsequently, various applications are briefly discussed for which hierarchies with different numbers of stages were developed.

\section{Abstract description}
The idea of a model hierarchy in the context of pa\-ra\-me\-trized partial differential equations~(PDEs) was originally introduced in~\cite{haasdonk2023certified}.
Here we describe the concept in a general form that is applicable in a wide range of scenarios and for several types of models.
\par
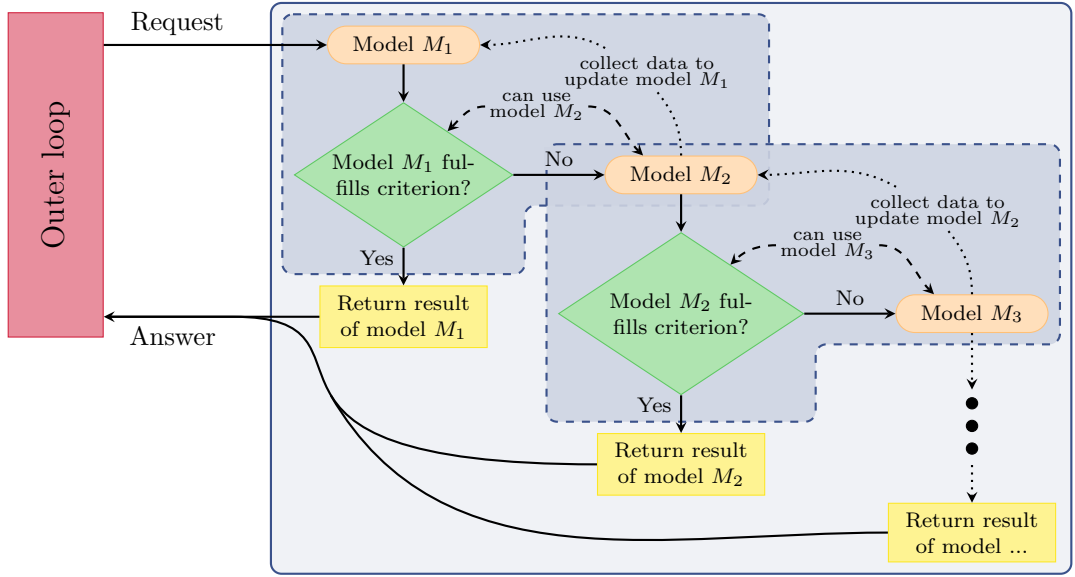
\begin{figure*}[htb]
	\centering
	\begin{tikzpicture}[align=center, node distance=0.5cm and 1.2cm, fill fraction/.style n args={2}{path picture={\fill[#1] ($(path picture bounding box.north west)!#2!(path picture bounding box.north east)$) rectangle (path picture bounding box.south east);}}]
	\node[model] (modelA) {Model~$M_1$};
	\node[error-estimator, below=of modelA] (criterionA) {Model~$M_1$ ful-\\fills criterion?};
	\node[rectangle, minimum width=1.25cm, minimum height=4.3cm, draw=matchingRed, fill=matchingRed!50, left=2.5cm of criterionA] (outerloop) {\rotatebox{90}{\large Outer loop}};
	\node[model, right=of criterionA] (modelB) {Model~$M_2$};
	\node[output, below=of criterionA] (outputA) {Return result\\of model~$M_1$};

	\draw[arrow] (outerloop.east |- modelA.west)  -- (modelA.west) node[above, pos=0.33] {Request};
	\draw[arrow] (modelA) -- (criterionA);
	\draw[arrow] (criterionA) -- (modelB) node[pos=0.5, yshift=6pt, inner sep=1.5pt, fill fraction={colorbg!25}{0.1}]{\small No};
	\draw[arrow] (criterionA) -- (outputA) node[pos=0.4, xshift=-10pt, yshift=2pt]{\small Yes};

	\node[inner sep=1.5pt] (modelB-data) at ($(modelB.north) + (-0.45cm, 1.1cm)$) {\scriptsize collect data to\\[-5pt]\scriptsize update model~$M_1$};
	\draw[dotted, thick] (modelB.north) to[out=90, in=310] ($(modelB-data.south) + (0.1cm, 0.06cm)$);
	\draw[arrow, dotted, thick] ($(modelB-data.north) + (-0.55cm, 0.05cm)$) to[out=150, in=0] (modelA.east);

	\node[inner sep=1.5pt] (modelB-usage) at ($(modelB.north) + (-1.9cm, 0.65cm)$) {\scriptsize can use\\[-5pt]\scriptsize model~$M_2$};
	\draw[arrow, dashed, thick] (modelB-usage.east) to[out=0, in=135] ($(modelB.north) + (-0.5cm, 0)$);
	\draw[arrow, dashed, thick] (modelB-usage.west) to[out=180, in=45] (criterionA);

	\node[error-estimator, below=of modelB, inner sep=3pt] (criterionB) {Model~$M_2$ ful-\\fills criterion?};
	\node[model, right=of criterionB] (modelC) {Model~$M_3$};
	\node[output, below=of criterionB] (outputB) {Return result\\of model~$M_2$};

	\draw[arrow] (modelB) -- (criterionB);
	\draw[arrow] (criterionB) -- (modelC) node[pos=0.5, yshift=6pt, inner sep=1.5pt]{\small No};
	\draw[arrow] (criterionB) -- (outputB) node[pos=0.4, xshift=-10pt, yshift=2pt]{\small Yes};

	\node[inner sep=1.5pt] (modelC-data) at ($(modelC.north) + (-0.45cm, 1.1cm)$) {\scriptsize collect data to\\[-5pt]\scriptsize update model~$M_2$};
	\draw[dotted, thick] (modelC.north) to[out=90, in=310] ($(modelC-data.south) + (0.1cm, 0.06cm)$);
	\draw[arrow, dotted, thick] ($(modelC-data.north) + (-0.55cm, 0.05cm)$) to[out=150, in=0] (modelB.east);

	\node[inner sep=1.5pt] (modelC-usage) at ($(modelC.north) + (-1.9cm, 0.65cm)$) {\scriptsize can use\\[-5pt]\scriptsize model~$M_3$};
	\draw[arrow, dashed, thick] (modelC-usage.east) to[out=0, in=135] ($(modelC.north) + (-0.5cm, 0)$);
	\draw[arrow, dashed, thick] (modelC-usage.west) to[out=180, in=45] (criterionB);
	\draw[arrow, dotted, thick] (modelC.south) -- ($(modelC.south)+(0,-0.75cm)$);

	\node[below=0.85cm of modelC, draw=black, fill=black, inner sep=1.5pt, circle] {};
	\node[below=1.15cm of modelC, draw=black, fill=black, inner sep=1.5pt, circle] {};
	\node[below=1.45cm of modelC, draw=black, fill=black, inner sep=1.5pt, circle] {};

	\node[output, below=2.25cm of modelC] (outputC) {Return result\\of model~...};
	\draw[arrow, dotted, thick] ($(outputC.north)+(0,0.5cm)$) -- (outputC.north);

	\draw[arrow] (outputA.west)  -- (outerloop.east |- outputA.west) node[below, pos=0.68] {Answer};
	\coordinate (temp) at ($(outputA.west)!0.4!(outerloop.east |- outputA.west)$);
	\draw[thick] (outputB) .. controls ($(outputB -| temp)+(-0.5cm,0)$) and ($(outputA)+(0,-0.cm)$) .. (temp) node[pos=0.49] (temp2) {} -- ($(outerloop.east |- outputA.west)+(0.075cm,0)$);
	\draw[thick] (outputC) to[out=180, in=-60] ($(temp2)+(0.5pt,-1pt)$);

	\begin{pgfonlayer}{background}
		\path[fill=colorbg!25, rounded corners, draw=colorbg, thick, fill opacity=0.3] ($(criterionA.west |- modelA.north)+(-2*\distbg,2*\distbg)$) -- ($(modelC.east |- modelA.north)+(2*\distbg,2*\distbg)$) -- ($(modelC.east |- outputC.south)+(2*\distbg,-\distbg)$) -- ($(outputC.south -| criterionA.west)+(-2*\distbg,-\distbg)$) -- cycle;
	\end{pgfonlayer}%

	\begin{pgfonlayer}{background}
		\path[fill=colorbg!25, rounded corners, draw=colorbg, dashed, thick, fill opacity=0.9] ($(criterionA.west |- modelA.north)+(-\distbg,\distbg)$) -- ($(modelA.north -| modelB.east)+(\distbg,\distbg)$) -- ($(modelB.south east)+(\distbg,-\distbg)$) -- ($(criterionA.east |- modelB.south)+(\distbg,-\distbg)$) -- ($(outputA.north -| criterionA.east)+(\distbg,\distbg)$) -- ($(outputA.north -| criterionA.west)+(-\distbg,\distbg)$) -- cycle;
	\end{pgfonlayer}%

	\begin{pgfonlayer}{background}
		\path[fill=colorbg!25, rounded corners, draw=colorbg, dashed, thick, fill opacity=0.9] ($(criterionB.west |- modelB.north)+(-\distbg,\distbg)$) -- ($(modelB.north -| modelC.east)+(\distbg,\distbg)$) -- ($(modelC.south east)+(\distbg,-\distbg)$) -- ($(criterionB.east |- modelC.south)+(\distbg,-\distbg)$) -- ($(outputB.north -| criterionB.east)+(\distbg,\distbg)$) -- ($(outputB.north -| criterionB.west)+(-\distbg,\distbg)$) -- cycle;
	\end{pgfonlayer}%
\end{tikzpicture}%
	\caption{Visualization of an abstract model hierarchy applied within an outer loop that sends requests to the hierarchy.}
	\label{fig:model-hierarchy-in-outer-loop}
\end{figure*}

In our abstract description we consider a solution operator~$\solutionMap\colon\inputSpace\to\solutionSpace$ that maps from an admissible input space~$\inputSpace$ to a possibly infinite dimensional solution space~$\solutionSpace$, where usually we know that~$\solutionMap$ exists, but it might not be accessible.
A typical example would be the solution operator of a parameterized~PDE where~$\inputSpace$ corresponds to the parameter set.
Furthermore, we assume that we are given a hierarchy of approximate models~$M_1, M_2,\dots$ that approximate the map~$\solutionMap$, where two successive models~$M_l$ and~$M_{l+1}$ in the hierarchy satisfy the following multi-fidelity assumptions:
\begin{itemize}
	\item  $\complexity{M_l} <  \complexity{M_{l+1}}$, where~$\complexity{M_l}$ denotes the model complexity as a measure for the runtime.
	\item  $\errorMeasure{M_l(\mu),\mu} \geq \errorMeasure{M_{l+1}(\mu),\mu}$, where~$\errorMeasure{M_l(\mu),\mu}$ denotes an error measure w.r.t.~$\solutionMap(\mu)$ for~$\mu\in\inputSpace$.
	\item  Model~$M_l$ can be improved by means of information from model~$M_{l+1}$.
\end{itemize}
Assume now that a request~$\mu \in\inputSpace$ in an outer multi-query loop needs to be processed.
The request is first passed to~model~$M_1$ which produces a result~$M_1(\mu)$ for the request.
This result is evaluated using the error measure, i.e.~it is verified whether~$\errorMeasure{M_1 (\mu),\mu} \leq {\rm TOL}$ is satisfied.
In order to check the criterion it might be necessary to also retrieve additional information from~model~$M_{2}$.
In general, if~model~$M_l$ fulfills the criterion, the result of~model~$M_l$ is returned to the outer loop.
If the criterion is not met, the request is passed to~model~$M_{l+1}$ which is assumed to be more accurate and is therefore more likely to fulfill the prescribed criterion.
Model~$M_{l+1}$ now proceeds similar to~model~$M_l$, i.e.~the request is processed resulting in an answer of~model~$M_{l+1}$.
When evaluating~model~$M_{l+1}$, data is collected that can be used, according to the third assumption from above, to improve~model~$M_l$.
Hence, model~$M_l$ is constructed and enhanced in an adaptive manner.
The result of~$M_{l+1}$ might now be passed on, depending on the structure of the remaining parts of the hierarchy.
Due to the involved check of the accuracy criterion for all results, the output of the model hierarchy is certified.
The overall hierarchical structure of the multi-fidelity algorithm is shown in~Fig.~\ref{fig:model-hierarchy-in-outer-loop} when applied in an outer loop for a hierarchy consisting of multiple stages.
For the algorithm performing the outer loop, the hierarchy behaves like a single model that returns a certain result of guaranteed accuracy.
All the internal model selection and adaptation is invisible from the outside.

\section{Applications}
In the following paragraphs we discuss applications where the concept of an adaptive model hierarchy was utilized successfully to speed up different computational tasks.
\paragraph*{PDE-constrained optimization.}
In~\cite{keil2022adaptive} we introduced a two-stage hierarchy consisting of a full order model and a machine learning surrogate for a PDE-constrained optimization problem that occurs in enhanced oil recovery.
The machine learning surrogate approximates the objective function and is based on training data gathered when evaluating the full order model (involving the costly simulation of a three-phase flow in a porous medium).
From the point of view of the model hierarchy as shown in~Fig.~\ref{fig:model-hierarchy-in-outer-loop}, the machine learning surrogate corresponds to~model~$M_1$ and during its evaluation an optimization problem for the approximate objective function is solved.
The accuracy of the result of this inner optimization loop is then evaluated by computing an approximate gradient using the full model~$M_2$.

\paragraph*{Parametrized parabolic PDEs.}
As a second example, we considered in~\cite{haasdonk2023certified} parametrized parabolic~PDEs where the hierarchy consists of a full order model, a reduced basis reduced order model and a machine learning model.
The latter model is based on the approach of learning the reduced coefficients with respect to a reduced basis as introduced in~\cite{hesthaven2018nonintrusive}.
The reduced basis is computed using evaluations of the full order model whereas the machine learning surrogate is trained based on solutions of the reduced basis model.
Accuracy of the reduced basis and the machine learning model is verified by means of an a posteriori error estimator for reduced models of parabolic problems.
Since the machine learning surrogate uses the same reduced space as the reduced basis model, the a posteriori error estimator is applicable also to the machine learning approximation.
Hence, a close connection between the two surrogate models facilitates their interaction in the hierarchy in this case.
The full order model here serves as reference and is therefore assumed to be arbitrarily accurate.
Hence, no accuracy check of the full order solution is performed.

\paragraph*{Parametrized optimal control problems.}
A three-stage adaptive model hierarchy for linear-quadratic optimal control problems with parameter-dependent system components was developed in~\cite{kleikamp2024application}.
The general structure is comparable to the one for parabolic~PDEs.
In particular, the three involved models and their interaction are similar and an a posteriori error estimator is used to certify the results obtained by the reduced models.
The special structure of the considered optimal control problems allows to identify solutions to the associated optimality system by the optimal adjoint at final time.
The reduced basis model thus builds on an approximation of the set of optimal final time adjoints by linear subspaces.
As before, the machine learning surrogate makes use of the same reduced space which allows to reuse the a posteriori error estimate of the reduced basis model.

\section{Conclusion}
The introduced concept of adaptive model hierarchies provides a possibility to combine different models of varying complexity within a joint hierarchy that can be evaluated efficiently.
As discussed in the last section, model hierarchies are applicable in different contexts and make use of the advantages of all involved models.

\bibliography{ifacconf}             
\end{document}